\documentclass[12pt,reqno]{article}

\usepackage[usenames]{color}
\usepackage{amssymb}
\usepackage{graphicx}
\usepackage{amscd}

\usepackage[colorlinks=true,
linkcolor=webgreen,
filecolor=webbrown,
citecolor=webgreen]{hyperref}

\definecolor{webgreen}{rgb}{0,.5,0}
\definecolor{webbrown}{rgb}{.6,0,0}

\usepackage{color}
\usepackage{fullpage}
\usepackage{float}

\usepackage{graphics,amsmath,amssymb}
\usepackage{amsthm}
\usepackage{amsfonts}
\usepackage{latexsym}

\begin{document}

\vspace*{2.1cm}

\theoremstyle{plain}
\newtheorem{theorem}{Theorem}
\newtheorem{corollary}[theorem]{Corollary}
\newtheorem{lemma}[theorem]{Lemma}
\newtheorem{proposition}[theorem]{Proposition}
\newtheorem{obs}[theorem]{Observation}
\newtheorem{claim}[theorem]{Claim}

\theoremstyle{definition}
\newtheorem{definition}[theorem]{Definition}
\newtheorem{example}[theorem]{Example}
\newtheorem{remark}[theorem]{Remark}
\newtheorem{conjecture}[theorem]{Conjecture}

\begin{center}


{\Large\bf Total coloring and efficient domination applications to non-Cayley non-Schreier vertex-transitive graphs} 

\vskip 1cm
\large
Italo J. Dejter

University of Puerto Rico

Rio Piedras, PR 00936-8377

\href{mailto:italo.dejter@gmail.com}{\tt italo.dejter@gmail.com}
\end{center}

\begin{abstract} Let $0<k\in\mathbb{Z}$. Let the star 2-set transposition graph $ST^2_k$ be the $2(k-1)$-regular graph whose vertices are the $2k$-strings on $k$ symbols, each symbol repeated twice, with its edges given each by the transposition of the initial entry of one such $2k$-string with any entry that contains a different symbol than that of the initial entry. The pancake 2-set transposition graph $PC^2_k$ has the same vertex set of $ST^2_k$ and its edges involving each the maximal product of concentric disjoint transpositions in any prefix of an endvertex string, including the external transposition being that of an edge of $ST^2_k$. For $1<k\in\mathbb{Z}$, we show that $ST^2_k$  and $PC^2_k$, among other intermediate transposition graphs, have total colorings via $2k-1$ colors. They, in turn, yield efficient dominating sets, or E-sets, of the vertex sets of $ST^2_k$ and $PC^2_k$, and partitions into $2k-1$ such E-sets, generalizing Dejter-Serra work on E-sets in such graphs. \end{abstract}

\section{Efficient domination and total coloring of graphs}\label{s1}

Let $0<k\in\mathbb{Z}$. Given a finite graph $G=(V(G),E(G))$ and a subset $S\subseteq V(G)$, it is said that $S$ is an {\it efficient dominating set} (E-set) \cite{AK,BBS,gen,D73,D76,edig} or a {\it perfect code} \cite{Borges,D80}, if for each $v\in V(G)\setminus S$ there exists exactly one vertex $v^0$ in $S$ such that $v$ is adjacent to  $v^0$.   

Applications of E-sets occur in: {\bf(a)} the theory of error-correcting codes and {\bf(b)} establishing the existence of regular graphs for Network Theory by removing E-sets from their containing graphs.

A {\it total coloring} of a graph $G$ is an assignment of colors to the vertices and edges of $G$ such that no two incident or adjacent elements (vertices or edges)
are assigned the same color \cite{tc-as}. 
A total coloring of $G$ such that the vertices adjacent to each $v\in V(G)$ together with $v$ itself are assigned pairwise different colors will be said to be an {\it efficient coloring}. Such efficient coloring will be said to be {\it totally efficient} if $G$ is regular, the color set is $[k]=\{0,1,\ldots,k-1\}$
and each $v\in V(G)$ together with its neighbors are assigned all the colors in $[k]$.
The {\it total} (resp. {\it efficient}) {\it chromatic number} $\chi''(G)$ (resp. $\chi'''(G)$) of $G$ is defined as the least number of colors required by a total (resp. efficient) coloring of $G$. 

As for applications other than (a)-(b) above, note that: {\bf(c)} by removing the vertices of a fixed color, then again regular graphs for Network Theory are generated; {\bf(d)} by removing the edges of a fixed color, then copies of a non-bipartite biregular graph whose parts have vertices with degrees differing in a unit are determined, again applicable in Network Theory.

In Section~\ref{s3}, we show that the graphs of a family of graphs $G=ST^2_k$, ($0<k\in\mathbb{Z}$), introduced in Section~\ref{s2}, satisfy the conditions of the following theorem. We conjecture that those conditions are only satisfied by such graphs $G=ST^2_k$, and not any other graphs.

\begin{theorem}\label{chi} {\bf(I)} Let $3<h\in 2\mathbb{Z}$. Let $G$ be a connected $(h-2)$-regular graph with a totally efficient coloring via color set $[h]\setminus\{0\}=\{1,\ldots,h-1\}$. Then, there is a partition of $V$ into $h-1$ subsets $W_1,\ldots,W_{h-1}$, where $W_i$ is formed by those vertices of $G$ having color $i$, for each $i\in[h]\setminus\{0\}$. In such a case, $\chi'''(G)=h-1$. Moreover, each $W_i$ is an E-set of $G$, for $i\in[h]\setminus\{0\}$. {\bf(II)} Let $4<h\in 2\mathbb{Z}$. Then, $G\setminus W_i$ is a connected $(h-3)$-regular subgraph that still has efficient chromatic number $h-1$, i.e. $\chi'''(G\setminus W_i)=h-1$, even though it has only a (non-total) efficient coloring. Letting $E_i$ be the set of edges with color $i$ in $G\setminus W_i$, then:\begin{enumerate} \item $G\setminus W_i\setminus E_i$ is the disjoint union of copies of regular subgraphs of degree $h-4$ with efficient colorings by $h-3$ colors obtained from $[h]\setminus\{0,i\}$ by removing the edges of a color $j\ne i$;
\item $G\setminus E_i$ is a non-bipartite $(h-2,h-3)$-biregular graph.
\end{enumerate} 
\end{theorem}

\begin{proof} We use the inequality $\chi''(G)\ge\Delta(G)+1$, where $\Delta(G)$ is the maximum degree of $G$ \cite{tc-as}. In our case, $\chi'''(G)=\chi''(G)=\Delta(G)+1$. Because of this,
a totally efficient coloring here provides a partition $W_1,\ldots,W_{h-2}$ as claimed in item (I). By definition of totally efficient coloring, each $W_i$ is an E-set. For item (II), deleting $W_i$ from $G$ removes also all the edges incident to the vertices of $W_i$, so $G\setminus W_i$ still has an efficient coloring which is not totally efficient since there is an edge color lacking incidence to each particular vertex of $G\setminus W_i$.
To establish item (II)1, note that removal of $E_i$ from $G\setminus W_i$ for $h>4$, leaves us with the graph induced by the edges of all colors other than color $i$, which necessarily disconnects $G\setminus W_i$, again because of the definition of totally efficient coloring. 
To  establish item (II)2, the removal of the edges with color $i$ leaves their endvertices with degree $h-3$ and forming a vertex subset of the resulting $G\setminus E_i$, while the remaining vertices have color $i$, degree $h-2$ and form a stable vertex set.
This completes the proof of the theorem.
 All of this can be verified without loss of generality via the proof of Theorem~\ref{t2}, for $h=2k$. 
\end{proof}

Let $\ell\in\{0,1\}$. In Section~\ref{Cs}, we generalize via $\ell$-set permutations, (see Section~\ref{s2}), the result of \cite{D73} that the star transposition graphs form a {\it dense segmental neighborly E-chain}. In Section~\ref{panqueque}, we generalize star transposition graphs to pancake transposition graphs and related intermediate graphs \cite{D73}, leading to an adequate version of dense neighborly E-chain \cite{D73}, with obstructions preventing any convenient version of segmental E-chain \cite{D73}. 

\section{Families of multiset transposition graphs with E-sets}\label{s2}

\begin{figure}[htp]
\includegraphics[scale=0.86]{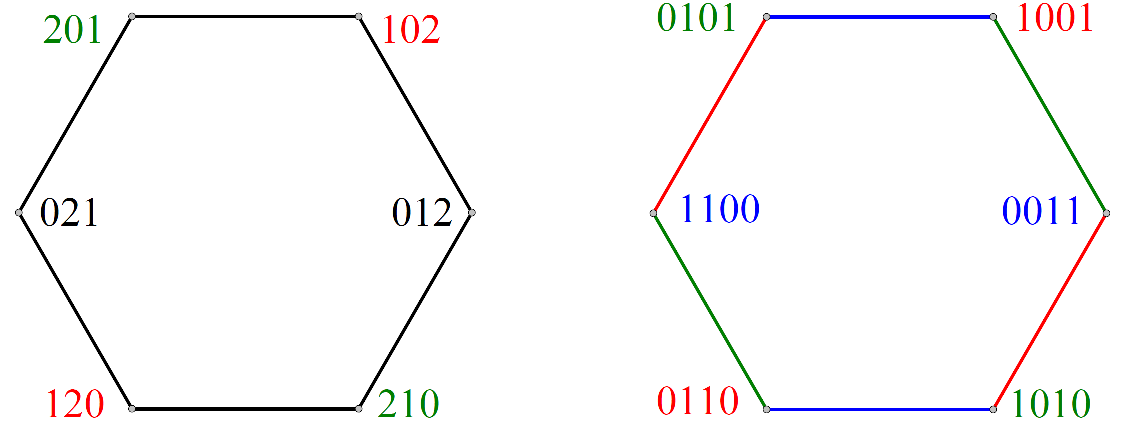}
\caption{The 6-cycles $ST^1_3=ST_3$ and $ST^2_2$.}
\label{fig1}
\end{figure}

Let $0<\ell\in\mathbb{Z}$ and let $1<k\in\mathbb{Z}$. We say that a string over the alphabet $[k]$ that contains exactly $\ell$ occurrences of $i$, for each $i\in[k]$, is an $\ell$-{\it set permutation}.
In denoting specific $\ell$-set permutations, commas and brackets are often omitted. 

Let $V^\ell_k$ be the set of all $\ell$-set permutations of length $k\ell$. Let the {\it star $\ell$-set transposition graph} $ST^\ell_k$ be the graph on vertex set $V^\ell_k$ with an edge between each two vertices $v=v_0v_1\cdots v_{k\ell-1}$ and $w=w_0w_1\cdots w_{k\ell-1}$ that differ in a {\it star transposition}, i.e. by swapping the first entry $v_0$ of $v=v_0v_1\cdots v_{k\ell-1}\in V^\ell_k$ with any entry $v_j$ ($j\in[k\ell]\setminus\{0\}$) whose value differs from that of  $v_0$ (so $v_j\ne v_0$), thus obtaining either $w=w_0\cdots w_j\cdots w_{k\ell-1}=v_j\cdots v_0\cdots w_{k\ell-1}$ or $w=w_0\cdots w_{k\ell-1}=v_{k\ell-1}\cdots v_0$. In other words,
each edge of $ST^\ell_k$ is given by the transposition of the initial entry of an endvertex string with an entry that contains a different symbol than that of the initial entry. The graphs $ST_k^\ell$ are a particular case of the graphs treated in \cite{faltaba} in a context of determination of Hamilton cycles.

It is known that all $k$-permutations, (that is all 1-set permutations of length $k$), form the {\it symmetric group}, denoted $Sym_k$, under composition of
$k$-permutations, each $k$-permutation $v_0v_1\cdots v_{k-1}$ taken as a bijection from the {\it identity} $k$-permutation $01\cdots(k-1)$ onto $v_0v_1\cdot v_{k-1}$ itself. 
A graph $ST^1_k$ with $k>1$ (which excludes $ST^1_1$) is the Cayley graph of $Sym_k$ with respect to the set of transpositions $\{(0\;i); i\in[k]\setminus\{0\}\}$.
Such a graph is denoted $ST_k$ in \cite{AK, D73}, where is proven its vertex set admits a partition into $k$ E-sets, exemplified on the left of Figure~\ref{fig1}
for $ST^1_3=ST_3$, with the vertex parts of the partition differentially colored in black, red and green, for respective first entries 0, 1 and 2. Figure 1 of \cite{D73} shows a similar example for $ST^1_4=ST_4$. Also, the graph $ST^2_k$ is vertex transitive, but is not a Cayley graph; see Subsection~\ref{den}, below.

\section{E-sets of star 2-set transposition graphs}\label{s3}

Let $i\in[2k]\setminus\{0\}=\{1,\ldots,2k-1\}$.
Let $\Sigma_i^k$ be the set of vertices $v_0v_1\cdots v_{k\ell-1}$ of $ST^\ell_k$ such that $v_0=v_i$, ($i=1,\ldots,2k-1$). 
Let $E_i^k$ be the set of edges having color $i$ in $G\setminus\Sigma_i^k$.
We will show that $\Sigma_i^k$ is an E-set of $ST^2_k$. Clearly, no edge of $E_i^k$ is incident to the vertices of $\Sigma_i^k$.

\begin{theorem}\label{t2} Let $k>1$. {\bf(i)}
The graph $ST^2_k$ has $\frac{(2k)!}{2^k}$ vertices and regular degree $2(k-1)$. {\bf(ii)} Let $i\in[2k]\setminus\{0\}=\{1,\ldots,2k-1\}$ and
let $\Sigma_i^k$ be the set of vertices $v_0v_1\dots v_{2k-1}$ of $ST^2_k$ such that $v_0=v_i$. Then, $V(ST^2_k)$ admits a partition into $2k-1$ E-sets $\Sigma_i^k$, ($i\in[2k]\setminus\{0\}$). {\bf(iii)} 
Let $k>2$, let $j\in[2k]\setminus\{0\}$ and let
$E_j^k$ be the set of all edges of color $j$. Then, $ST^2_k\setminus\Sigma_i^k\setminus E_i^k$ is the disjoint union of $k2^{k-1}$ copies of $ST^2_{k-1}$. 
\end{theorem}

\begin{proof} Let $i=2k-1$ and let $j\in[2k]$. Then, 
each vertex $v=v_0v_1\cdots v_{2k-3}v_{2k-2}v_{2k-1}=0v_1\cdots v_{2k-3}j0$ is the neighbor of vertex $w=jv_1\cdots v_{2k-3}00$ via an edge of color $k-1$.
But $v\in\Sigma_i^k=\Sigma_{2k-1}^k$. Being $w$ at distance 1 from $\Sigma_{2k-1}^k$, then $w$ is in the {\it open neighborhood} $N(\Sigma_i^k)$ \cite{D73} of $\Sigma_{2k-1}^k$ in $ST^2_k$, so $w\in N(\Sigma_i^k)=N(\Sigma_{2k-1}^k)\subseteq ST^2_k\setminus\Sigma_i^k\setminus E_i^k=ST^2_k\setminus\Sigma_{2k-1}^k\setminus E_{2k-1}^k$. In fact, $N(\Sigma_i^k)=N(\Sigma_{2k-1}^k)$ is a connected component of $ST^2_k\setminus\Sigma_i^k\setminus E_i^k=ST^2_k\setminus \Sigma_{2k-1}^k\setminus E_{2k-1}^k$. A similar conclusion holds for each other open neighborhoods $N(\Sigma_i^k)$, ($1\le i<2k-1$).
\end{proof}

\begin{remark}\label{occur} The total coloring of $ST^2_k$ will be referred to as its {\it color structure}. 
The $k2^{k-1}$ copies of $ST^2_{k-1}$ in $ST^2_k$ whose disjoint union is $ST^2_k\setminus\Sigma_i^k\setminus E_i^k$ inherit each a color structure that generalizes that of Examples~\ref{olvid}-\ref{olvido}, below, and is similar to the color structure of $ST^2_{k-1}$.
\end{remark}

\begin{example}\label{olvid}
The graph $ST^2_2$ has the totally efficient coloring depicted on the right of Figure~\ref{fig1}, where $\Sigma_1^2=\{0011,1100\}$ is color blue, as is $E_1^2=\{(0101,1001),(0110,1010)\}$; $\Sigma_2^2=\{0101,1010\}$ is color green, as is $E_2^2=\{(0110,1100),(0011,1001)\}$; $\Sigma_3^2=\{0110,1001\}$ is color red, as is $E_3^2=\{(0011,1010),(0101,1100)\}$.\end{example}

\begin{example}\label{olvido}
The graph $ST^2_3$ has the E-set $\Sigma_5^3$ with 18 vertices denoted as in display (\ref{18}):  
\begin{eqnarray}\label{18}\begin{array}{cccccc}
A=011220,&\underline{A}=022110,&B=012210,&\underline{B}=021120,&C=012120,&\underline{C}=021210,\\
D=122001,&\underline{D}=100221,&E=120021,&\underline{E}=102201,&F=120201,&\underline{F}=102021,\\
G=200112,&\underline{G}=211002,&H=201102,&\underline{H}=210012,&J=201012,&\underline{J}=210102.\\
\end{array}\end{eqnarray}

\begin{figure}[htp]
\includegraphics[scale=0.88]{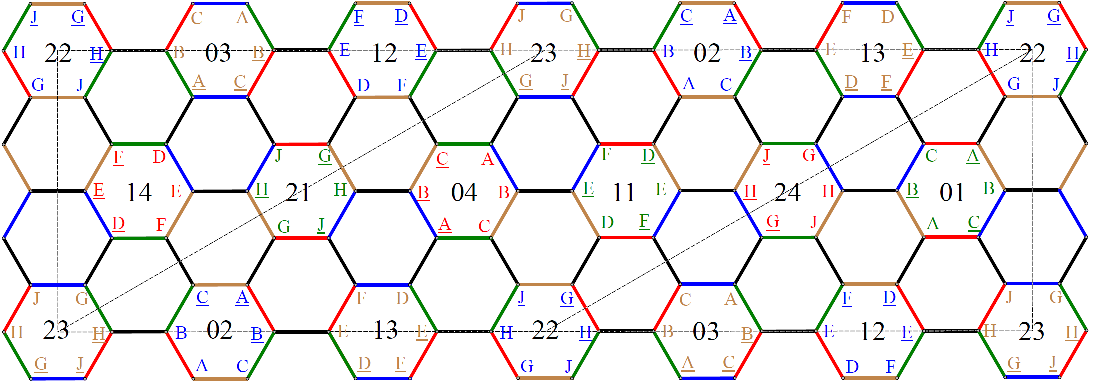}
\caption{A fundamental region of a lattice suggests a rhomboidal torus cutout of $ST^2_3$.}
\label{fig3}
\end{figure}

A planar interconnected disposition of the 6-cycles of the subgraph $ST^2_3\setminus \Sigma_5^3$ of $ST^2_3$ is shown in Figure~\ref{fig3}. The edges of such 6-cycles are alternatively colored with 2 or 3 colors of the color form $(ababab)$ or $(abcabc)$ respectively, where $\{a,b,c\}\subseteq\{1,2,3,4\}$  is a subset of colors provided by the respective positions 1,2,3,4 of the $6$-tuples taken as the vertices of $ST^2_3$. 

The tessellation suggested in Figure~\ref{fig3} can be extended to the whole plane as an unfolding of the fundamental region delimited by the shown dash-border rectangle -- call it $R$. This $R$ appears partitioned via dashed segments into two right triangles and a rhomboid in between. By transporting the left right triangle -- call it $T_l\subset R$ -- to a new position $T'_l$ to the right so that the vertical side of $T'_l$ coincides with the right side of $R$, a rhomboid $R'$ is obtained. Identification of the tilted sides of $R'$ and of its horizontal sides allows to view a toroidal embedding of $ST^2_3\setminus \Sigma_5^3$. 

\begin{figure}[htp]
\includegraphics[scale=0.87]{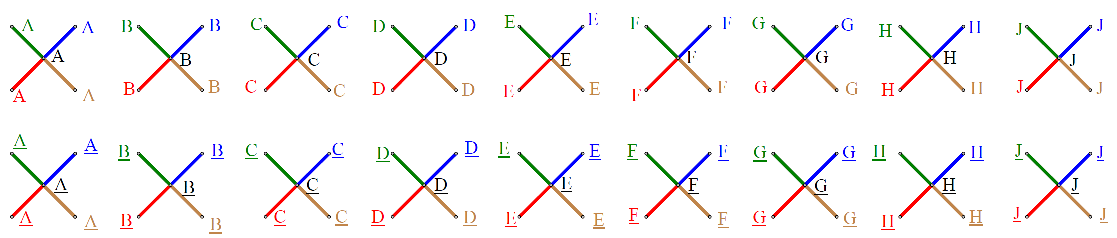}
\caption{The eighteen stars $K_{1,4}$ in $ST^2_3$ centered at the vertices of the E-set $\Sigma_5^3$.}
\label{fig4}
\end{figure}

Edge colors in Figure~\ref{fig3} are numbered as follows (indicating corresponding subsequent positions in the 6-tuples representing the vertices of $ST^2_3$): 
\begin{eqnarray}\label{AJ}1=green,\; 2=blue,\; 3=hazel,\; 4=red,\; 5=black.\end{eqnarray}
In Figure~\ref{fig3}, the 3-colored 6-cycles are exactly those containing in their interiors (next to their corresponding denoting vertices) the (possibly underlined) capital letters of display (\ref{18}), but each such letter colored as indicated in display (\ref{AJ}). Each such number color $a\in\{1,2,3,4\}$ as in display (\ref{AJ}) of a symbol $X\in\{A,\ldots,J,\underline{A},\ldots,\underline{J}\}$ in Figure~\ref{fig3} indicates the existence of an (absent) $a$-colored edge between $V^2_3\setminus\Sigma_5^3$ and $\Sigma_5^3$ in $ST^2_3$. Figure~\ref{fig4} shows each such edge in exactly one copy $\Upsilon$ of $K_{1,4}$ with its endvertex in $\Sigma_5^3$ represented by $X$ (in black) and its other endvertex being the sole element of  $\Upsilon\cap V^2_3\setminus\Sigma_5^3$, namely the $a$-colored $X$, that we denote as $X^a$ in Table~\ref{dibujito}. In fact, Table~\ref{dibujito} reproduces the data of Figure~\ref{fig3} in a likewise disposition, with the vertex notation $X^a$ instead of the $a$-colored $X$ notation of Figure~\ref{fig3}. In Table~\ref{dibujito}, edges are represented by their numeric symbols (display (\ref{AJ})) and appear interspersed with the symbols $X^a$ in representing the 3-colored 6-cycles, while 2-colored 6-cycles are represented by the disposition of their numeric symbols. Note in Figure~\ref{fig3} that each 3-colored 6-cycle is bordered by six 2-colored 6-cycles via edges colored in $\{1,2,3,4,\}$, while each 2-colored 6-cycle, call it $\Theta$, is bordered by three 3-colored 6-cycles (via edges in one fixed color of $\{1,2,3,4\}$) alternated with three 2-colored 6-cycles via an edge matching bordering $\Theta$ and whose color is 1.

\begin{table}[htp]
$$\begin{array}{lllllllllllllllllll}
1\underline{J}^23\underline{G}^24&&4C^32A^31&&1\underline{F}^23\underline{D}^24&&4J^32G^31&&1\underline{C}^23\underline{A}^24&&4F^32D^31&&\\
H^2\hspace*{6mm}\underline{H}^2&5& B^3\hspace*{6mm}\underline{B}^3&5&E^2\hspace*{6mm}\underline{E}^2&5&H^3\hspace*{6mm}\underline{H}^3&5&B^2\hspace{6mm}\underline{B}^2&5&E^3\hspace*{8mm}\underline{E}^3&5&\\
4G^23J^21&&1\underline{A}^32\underline{C}^34&&4D^23F^21&&1\underline{G}^32\underline{J}^34&&4A^23C^21&&1\underline{D}^32\underline{F}^34&&\\
5\hspace*{12mm}5&&5\hspace*{12mm}5&&5\hspace*{12mm}5&&5\hspace*{12mm}5&&5\hspace*{12mm}5&&5\hspace*{12mm}5\\
\end{array}$$
$$\begin{array}{rrrrrrrrrrrrrrrrrrr}
&&&3\underline{F}^41D^42&&2J^14\underline{G}^13&&3\underline{C}^41A^42&&2F^14\underline{D}^13&&3\underline{J}^41G^42&&2C^14\underline{A}^13\\
&&5&\underline{E}^4\hspace*{6mm}E^4&5&\underline{H}^1\hspace*{6mm}H^1&5&\underline{B}^4\hspace*{6mm}B^4&5&\underline{E}^1\hspace*{6mm}E^1&5&\underline{H}^4\hspace{6mm}H^4&5&\underline{B}^1\hspace*{6mm}B^1\\
&&&2\underline{D}^41F^43&&3G^14\underline{J}^12&&2\underline{A}^41C^43&&3D^14\underline{F}^12&&2\underline{G}^41J^43&&3A^14\underline{C}^12\\
\end{array}$$
$$\begin{array}{lllllllllllllllllll}
5\hspace*{12mm}5&&5\hspace*{12mm}5&&5\hspace*{12mm}5&&5\hspace*{12mm}5&&5\hspace*{12mm}5&&5\hspace*{12mm}5\\
4J^32G^32&&1\underline{C}^23\underline{A}^24&&4F^32D^31&&1\underline{J}^23\underline{G}^24&&4C^32A^31&&1\underline{F}^23\underline{D}^24&&\\
H^3\hspace*{6mm}\underline{H}^3&5&B^2\hspace{6mm}\underline{B}^2&5&E^3\hspace*{8mm}\underline{E}^3&5&H^2\hspace*{6mm}\underline{H}^2&5& B^3\hspace*{6mm}\underline{B}^3&5&E^2\hspace*{6mm}\underline{E}^2&5&\\
1\underline{G}^32\underline{J}^34&&4A^23C^21&&1\underline{D}^32\underline{F}^34&&4G^23J^21&&1\underline{A}^32\underline{C}^34&&4D^23F^21&&\\
\end{array}$$
\caption{Notational disposition of elements of $ST^2_3$ in Figure~\ref{fig3}.}
\label{dibujito}\end{table}

Table~\ref{los18} represents the twelve 3-colored 6-cycles, as follows. The six centers $X\in\{A,\ldots,J,$ $\underline{A},\ldots,\underline{J}\}$ of copies of $K_{1,4}$ involved with one such 3-colored 6-cycle, call it $\Phi$, are represented by 6-tuples that are expressed in Table~\ref{los18} in a 6-row section of a column whose heading is $\Sigma_5^3$. To the immediate right of each such 6-row section, another 6-row section of 6-tuples expresses the corresponding neighbors $X^b$, for a fixed color $b\in\{1,2,3,4\}$, via $b$-colored edges. 
Such neighbors $X^b$ conform $V(\Phi)$ and induce $\Phi$. In fact, Table~\ref{los18} contains the twelve instances of such representations. 

Notice that the vertices in display (\ref{18}) are of the form $ia_1a_2a_3a_4i$.
Centered inside each 3-colored 6-cycle $\Phi$ in Figure~\ref{fig3}, a pair $(i,b)$ of digits (written as $ib$) indicates the fixed double entry $i\in\{0,1,2\}$ of the vertices $ia_1a_2a_3a_4i$ of $\Sigma_5^3$ in $\Phi$ and the fixed color $b$ their representing symbols have in the figure. 

To facilitate viewing the edge colors along each $\Phi$, 
the second row in Table~\ref{los18} shows the 6-tuple $x$ of subsequent positions (or colors), $012345$, of the 6-tuples representing each $X$ and $X^b$. In each such $x$ under the heading $\Sigma_5^3$, the entry $b\in\{1,2,3,4\}$ of the corresponding $X^b$ is underlined, while under each subsequent heading $X^b$, the other three entries in $\{1,2,3,4\}$ are underlined to indicate the  
entries successively transposed with the initial entry in the subsequent vertically disposed 6-tuples of each particular $\Phi$. 

Observe the difference between 3-colored 6-cycles appearing here and 2-colored 6-cycles in that the former are created by transpositions not involving the initial entry while the latter do involve transpositions with the initial entry.

In Figure~\ref{fig3}, deletion of the edges colored 1 from $ST^2_3\setminus\Sigma_5^3$ leaves a subgraph with twelve components, each being a 3-colored 6-cycle. Note that $E(ST^2_3)$ has a 1-factorization into five 1-factors $E_1^3,E_2^3,E_3^3,E_4^3,E_5^3$, each $E_i^3$ composed by those edges colored $i$, ($i\in[6]\setminus\{0\}$). Moreover, $ST^2_3\setminus\Sigma_5^3\setminus E_5^3$ is the union of the twelve 3-colored 6-cycles in Table~\ref{los18}. 
\end{example}

\begin{corollary}\label{cor} Let $k>2$. Then:
\begin{enumerate}
\item $ST^2_k$ has $\frac{2k!}{2^k}$ vertices having $\frac{2k!}{2^k(2k-1)}$ vertices in each color $1,2,\ldots,2k-1$;
\item $ST^2_k$ has $\frac{2k!}{2^k}\times(k-1)$ edges;
\item color $k\ell-1$ provides exactly $\frac{2k!}{2^k(2k-1)}=y$ vertices forming an E-set  $\Sigma_{2k-1}^k$ of $ST^2_k$;
\item the $y$ resulting dominating copies of $K_{1,2k-2}$ have a total of $y\times(2k-2)=z$ edges;
\item there are exactly $\frac{2k!}{2^k}\times(k-1)-z=h$ edges in $ST_{2k-1}^k$ not counted in item 4;
\item the $h$ edges in item 5. contain $\frac{h}{2k-1}$ edges in each color $1,2,\ldots,2k-1$;
\item so they contain $h-\frac{h}{2k-1}$ edges in colors $\ne 2k-1$, (namely, $1,2,\ldots,2k-2$);
\item there are $\frac{2k!}{2^k}-y$ vertices in $ST^2_k\setminus\Sigma_{2k-1}^k$ dominated by $\Sigma_{2k-1}^k$;
\item the $\frac{2k!}{2^k}-y$ vertices in item 8. appear in $k\times(2k-2)$ copies of $ST^2_{k-1}$;
\item there are $\frac{h}{(2k-1)^2k}$ edges in each copy of $ST_{2k-1}^k$ in $ST^2_k\setminus\Sigma_{2k-1}^k$.
\end{enumerate}
\end{corollary}

\begin{table}[htp]
$$\begin{array}{||c|c|c||c|c|c||c|c|c||c|c|c||c|c|c||c|c|c|c|c||}\hline\hline
X&\Sigma_5^3&X^1&X&\Sigma_5^3&X^2&X&\Sigma_5^3&X^3&X&\Sigma_5^3&X^4\\\hline\hline
x&0\underline{1}2345&01\underline{234}5&x&01\underline{2}345&0\underline{1}2\underline{34}5&x&012\underline{3}45&0\underline{12}3\underline{4}5&x&0123\underline{4}5&0\underline{123}45\\\hline
A  &  011220    &    101220  &  A  &  011220  &  110220  &  A  &   011220  &  211020  &  A  &   011220  &  211200\\
\underline{B}  &  021120    &    201120  &  B  &  012210  &  210210  & \underline{B}  &  021120  &  121020  &  B  &  012210  &  112200\\
C  &  012120    &   102120  & \underline{C}  &  021210  &  120210  & \underline{C}  &  021210  &  221010  &  C  &  012120  &  212100\\
\underline{A}  &  022110    &    202110  & \underline{A}   &  022110  &  220110  & \underline{A}  &  022110  &   122010  & \underline{A} &  022110  &  122100\\
B  &  012210    &   102210  & \underline{B}  &  021120  &  120120  &  B  &  012210  &   212010  & \underline{B}  &  021120  &  221100\\
\underline{C}  &  021210    &   201210  &  C  &  012120  &  210120  &  C  &  012120  &  112020  & \underline{C}  &  021210  &  121200\\\hline

D  &  122001    &    212001  &  D  &  122001  &  221001  &  D  &   122001  &  022101  &  D  &  122001  &  022011\\
\underline{E}  &  102201    &    012201  &  E  &  120021  &  021021  & \underline{E}  &   102201  &  202101  &  E  &   120021  &  220011\\
F  &  120201    &    210201  & \underline{F}  &  102021  &  201021  & \underline{F}  &   102021  &  002121  &  F  &   120201  &  020211\\
\underline{D}  &  100221    &    010221  & \underline{D}  &  100221  &  001221  & \underline{D}  &  100221  &  200121  & \underline{D}  &  100221  &  200211\\
E  &  120021    &    210021  & \underline{E}  &  102201  &  201201  &  E  &   120021  &  020121  & \underline{E}  &  102201  &  002211\\
\underline{F}  &  102021    &    012021  &  F  &  120201  &  021201  &  F  &   120201  &  220101  & \underline{F}  &   102021  &  202011\\\hline
\
G  &  200112    &    020112  &  G  &  200112  &  002112  &  G  &  200112  &  100212  &  G  &  200112  &  100122\\
\underline{H}  &  210012    &    120012  &  H  &  201102  &  102102  &  \underline{H} &  210012  &  010212  &  H  &  201102  &  001122\\
J  &  201012    &     021012  & \underline{J}  &  210102  &  012102  & \underline{J}  &   210102  &  110202  &  J  &   201012  &  101022\\
\underline{G}  &  211002    &    121002  & \underline{G}  &  211002  &  112002  & \underline{G}  &  211002  &  011202  & \underline{G}  &  211002  &  011022\\
H  &  201102    &    021102  & \underline{H}  &  210012  &  012012  &  H  &  201102  &  101202  & \underline{H}  &  210012  &  110022\\
\underline{J}  &  210102    &    120102  &  J  &  201012  &  102012  &  J  &    201012  &  001212  & \underline{J}  &  210102  &  010122\\
\hline\hline
\end{array}$$
\caption{The twelve 6-cycles whose vertices start with 00, 11 and 22}
\label{los18}\end{table}

\begin{proof}
The ten items of the corollary can be verified directly from the enumerative facts involved with the graphs $ST^2_k$.
\end{proof}

\begin{example}\label{2^3}
For $ST^2_3$, we have that:
\begin{enumerate}
\item $ST^2_3$ has $\frac{6!}{2^3}=90$ vertices containing $\frac{90}{5}=18$ vertices in each color $1,2,3,4,5$;
\item $ST^2_3$ has $90\times 4/2=180$ edges;
\item color 5 provides 18 vertices that form an E-set  $\Sigma_5^3$ of $ST^2_3$;
\item the 18 resulting dominating copies of $K_{1,4}$ in $ST^2_3$ have $18\times 4=72$ edges;
\item outside that, there are still $180-72=108$ edges;
\item they contain $\frac{108}{5}=36$ edges in each color $1,2,3,4,5$;
\item so they contain $108-36=72$ edges in colors $\ne 5$, (namely, $1,2,3,4$);
\item there are $90-18=72$ remaining vertices in $ST^2_3$, dominated by $\Sigma_5^3$;
\item they appear in $3\times 4=12$ copies of $ST^2_2$;
\item there are $\frac{72}{3\times 4}=\frac{72}{12}=6$ edges in each copy of $ST^2_2$ in $ST^2_3\setminus\Sigma_5^3$.
\end{enumerate}
\end{example}

\begin{example}\label{2^4}
For  $ST^2_4$, we have that:
\begin{enumerate}
\item $ST^2_4$ has $\frac{8!}{2^4}=2520$ vertices containing $\frac{2520}{7}=360$ vertices in each color $1,\ldots,7$;
\item $ST^2_4$ has $2520\times 6/2=7560$ edges;
\item color 7 provides 360 vertices that form an E-set  $\Sigma_7^4$ of $ST^2_4$;
\item the 360 resulting dominating copies of $K_{1,6}$ in $ST^2_4$ have $360\times 6=2160$ edges;
\item outside that, there are still $7560-2160=5400$ edges;
\item they contain $\frac{5400}{7}=1080$ edges in each color $1,2,3,4,5,6,7$;
\item the $h$ edges in item 6 have $5040-1080=4320$ edges in colors $\ne 7$, (namely, $1,\ldots,6$);
\item there are $2520-360=2160$ remaining vertices in $ST^2_4$, dominated by $\Sigma_7^4$;
\item they appear in $4\times 6=24$ copies of $ST^2_3$;
\item there are $\frac{4320}{4\times 6}=\frac{4320}{24}=180$ edges in each copy of $ST^2_3$ in $ST^2_4\setminus\Sigma_7^4$.
\end{enumerate}
\end{example}

\begin{example}
The 24 copies of $ST^2_3$ in $ST^2_4$, (item 5 of Example~\ref{2^4}), can be encoded as follows. We start by encoding the fundamental rectangle in Figure~\ref{fig3} by arranging the pairs $(i,b)=ib$ as follows, following the disposition in the figure:

\begin{eqnarray}\label{3x7}\begin{array}{ccccccccccccc}
22&&03&&12&&23&&02&&13&&22\\
&14&&21&&04&&11&&24&&01&\\
23&&02&&13&&22&&03&&12&&22\\
\end{array}\end{eqnarray}

By further encoding this disposition as $(012,1234)$, we now have that the 24 copies of $ST^2_3$ in $ST^2_4$ can be expressed as:
$$(123,123456), (013,123456), (023,123456), (012,123456).$$
\end{example}

A characterization of the twenty-four 2-colored 6-cycles of $ST^2_3\setminus\Sigma_1^3$ is also available from that of the twelve 3-colored 6-cycles in display (\ref{3x7}). Let us observe the triple $(0x_0,1y_1,2y_2)$ formed by the three pairs $0x_0$, $1x_1$, $2x_2$ denoting the three 3-colored 6-cycles that share each an edge $e$ with a given 2-colored 6-cycle $\Theta_e$. By shortening each such triple of pairs to the triple of colors $x_0x_1x_2$ and setting its missing color $x_3$ in $\{1,2,3,4\}$ as a subindex, with colors $i=5$  and  $x_3$ assigned alternatively to the edges of each $\Theta_e$,
we have now the disposition in display (\ref{5pisos}) which is similar to that of  Figure~\ref{los18}: 
\begin{eqnarray}\label{5pisos}\begin{array}{ccccccccccccc}
22&&03&&12&&23&&02&&13&&22\\
142_3&342_1&341_2&321_4&421_3&423_1&413_2&213_4&214_3&234_1&134_2&132_4&142_3\\
&14&&21&&04&&11&&24&&01&\\
143_2&243_1&241_3&231_4&431_2&432_1&412_3&312_4&314_2&324_1&124_3&123_4&143_2\\
23&&02&&13&&22&&03&&12&&23\\
\end{array}\end{eqnarray}

Again, this disposition is encoded as $(123,1234)$. 

\begin{theorem}\label{app} The graphs $ST^2_k$ satisfy the conditions of Theorem~\ref{chi}, so they also satisfy its conclusions. 
\end{theorem}

\begin{proof} Because of the previous discussion, we see that in the hypotheses of Theorem~\ref{chi}
it is enough to take $h=2k$, $G=ST^2_k$, $W_i=\Sigma_i^k$ and $E_i=E_i^k$.
\end{proof}

\section{Open Problems}\label{op} 
We conjectured that the graph $G$ in the statement of Theorem~\ref{chi} must necessarily coincide with some $ST^2_k$.
On the other hand, the twenty-four 2-colored 6-cycles of $ST^2_3\setminus\Sigma_5^3$ generalize to 2-colored 6-cycles in $ST^2_k\setminus\Sigma_{2k-1}^k$, for any $k>3$, by similarly alternating three black edges (meaning color $2k-1$) with three edges of a common color different from $2k-1$ in order to obtain one such 2-colored 6-cycle. Performing this to include all edges of $ST^2_k\setminus\Sigma_{2k-1}^k$, still   
we do not know how to generalize for $k>3$ what happens between the $k2^{k-1}$ copies of $ST^2_{k-1}$ in Theorem~\ref{t2} and the black edges (colored via $2k-1$). The determination of this particular matter is left as an open problem. 

As a hint to illuminate the problem, let us recall that
$ST^2_k$ has $\frac{(2k)!}{2^k}$ vertices and regular degree $2(k-1)$; then it has $\frac{(2k)!(k-1)}{2^k}$ edges and a total coloring via $2k-1$ colors. The number of vertices in $ST^2_k$ having a fixed color is $\frac{(2k)!}{2^k(2k-1)}$. The copies of stars $K_{1,2k-2}$ with centers on vertices of $ST^2_k$ having a fixed color contain a total of $\frac{(2k)!(2k-2)}{2^k(2k-1)}=\frac{2k)!(k-1)}{2^{k-1}(2k-1)}$ edges. The numbers of remaining vertices and edges, namely those of $ST^2_k\setminus\Sigma_{2k-1}^k$, are $\frac{(2k)!}{2^k}-\frac{(2k)!}{2^k(2k-1)}$ and $\frac{(2k)!(2k-1)}{2^k}-\frac{(2k)!(k-1)}{2^{k-1}(2k-1)}$, respectively.  The edges of $ST^2_k\setminus\Sigma_{2k-1}^k$ with a fixed color are divided into groups of three edges, each such group with alternate edges of a corresponding 2-colored 6-cycle, with the other three alternating edges in color $2k-1$. A conclusion here is that the number of 2-colored 6-cycles must be the third part of $\frac{(2k)!(2k-1)}{2^k}-\frac{(2k)!(k-1)}{2^{k-1}(2k-1)}$, which for $k=3$ equals 24, as can be counted for example via Figure~\ref{fig3}.

\section{Conclusions for star 2-set transposition graphs}\label{Cs}

Let us recall from \cite{D73} that: 
\begin{enumerate}  

\item a countable family of graphs $${\mathcal G}=\{\Gamma_1\subset\Gamma_2\subset\cdots\subset\Gamma_i\subset\Gamma_{i+1}\subset\cdots\}$$ is said to be an {\it E-chain} if every $\Gamma_i$ is an induced subgraph of $\Gamma_{i+1}$ and each $\Gamma_i$  has an E-set $C_i$;

\item for graphs $\Gamma$ and $\Gamma'$, a one-to-one graph homomorphism $\zeta:\Gamma\rightarrow\Gamma'$  
such that $\zeta(\Gamma)$ is an induced subgraph of $\Gamma'$ is said to be an {\it inclusive map}; 

\item for $i\ge 1$, let $\kappa_i$ be an inclusive map of $\Gamma_i$ into $\Gamma_{i+1}$; 
if $C_{i+1}=N(\kappa_i(V(\Gamma_i)))$, then the E-chain $\mathcal{G}$ is said to be a {\it neighborly} E-chain;

\item a particular case of E-chain $\mathcal{G}$ is the one in which $C_{i+1}$ has a partition into images $\zeta_i^{(j)} (C_i)$ of $C_i$ through respective inclusive maps $\zeta_i^{(j)}$, where $j$ varies on a suitable  finite indexing set. In such a case, the E-chain is said to be {\it segmental}. 

\end{enumerate}

The notion of neighborly E-chain in item 3 above is not suitable in our context of graphs $ST^2_k$ and their E-sets, that we denote $\Sigma_{2k-1}^k$ (instead of $C_i$ as in \cite{D73}), like 
$\Sigma_3^2$ and $\Sigma_5^3$ in Example~\ref{olvido}, with $\Sigma_5^3$ detailed both in display (\ref{18}) and Figures~\ref{fig3}-\ref{fig4}, and also in Tables~\ref{dibujito}-\ref{los18}. 
In this context, the graphs $ST^2_k$ form an E-chain
\begin{eqnarray}\label{G2}{\mathcal{ST}}(2)=\{ST^2_1\subset ST^2_2\subset\cdots\subset ST^2_k\subset ST^2_{k+1}\subset\cdots\},\end{eqnarray}
with each inclusion $ST^2_k\subset ST^2_{k+1}$ realized by a set of $k+1$ {\it neighborly}  maps
\begin{eqnarray}\label{z1}\kappa_k^j:ST^2_k\rightarrow ST^2_{k+1},\end{eqnarray} ($j\in[k+1]$), ({\it neighborly} meaning that the images $\kappa_k^j(ST^2_i)$ are pairwise disjoint in $ST^2_{k+1}$ and that
\begin{eqnarray}\label{z2}\Sigma_k^{k+1}=\cup_{j=1}^{k-1}N(\kappa_i^j(V^2_i))\end{eqnarray} as a disjoint union), these neighborly maps given by \begin{eqnarray}\label{z3}\kappa_k^j(a_0a_1\cdots a_{2k-2}a_{2k-1})=(a^j_0a^j_1\cdots a^j_{2k-2}a^j_{2k-1}jj),\end{eqnarray}
for $j\in[k+1]$, where 
\begin{eqnarray}\label{z4}a^k_i=a_i,\: a^{k+1}_i=a_i+1\mod(k+1),\; \ldots, a^{k+h}_i=a_i+h\mod(k+1), \ldots,\end{eqnarray} for $i=0.1,\ldots,2k-1$, the superindices $k+h$ of the entries $a^{k+h}_j$ taken mod $k+1$.

As an example, the last column of Table~\ref{los18} offers disjoint neighborly maps $\kappa_2^j$, for $j=0,1,2$, yielding respectively the following images of the 6-cycle that comprises $ST^2_2$: 
$$\begin{array}{c}
\kappa^2_2(1001,0011,1010,0110,1100,0101)=(100122,001122,101022,011022,110022,010122);\\
\kappa^0_2(1001,0011,1010,0110,1100,0101)=(211200,112200,212100,122100,221100,121200);\\
\kappa^1_2(1001,0011,1010,0110,1100,0101)=(022011,220011,020211,200211,002211,202011).\\
\end{array}$$

An E-chain as in display (\ref{G2}) where each inclusion $ST^2_k\subset ST^2_{k+1}$ is realized by $k+1$ neighborly maps $\kappa_k^j$, as defined in displays (\ref{z1}) to (\ref{z4}), is said to be a
{\it disjoint neighborly} E-chain. 

The notion of segmental E-chain can also be generalized to the case of the graphs $ST^2_k$, where in item 3 above we replace ``neighborly" by ``disjoint neighborly". In that case, the E-chain will be said to be {\it disjoint segmental}. It is clear by symmetry that the E-chain ${\mathcal{ST}}(2)$ of display (\ref{G2}) is disjoint segmental, as exemplified via Figures~\ref{fig3}-\ref{fig4} and the related Tables~\ref{dibujito}-\ref{los18}. 

If, for each $i\ge 1$, there exists an inclusive map $\zeta_i:\Gamma_i\rightarrow\Gamma_{i+1}$  such that $\zeta(C_i)\subset C_{i+1}$, then \cite{D73} calls the E-chain {\it inclusive} and observes that an inclusive neighborly E-chain has $\kappa_i\ne \zeta_i$, for every positive integer $i$. 

\subsection{Density}\label{den}

In addition, \cite{D73} calls an E-chain $\mathcal{G}$ {\it dense} if, for each $n\ge 1$, one has $|V(\Gamma_n)|=(n+1)!$ and $|C_n|=n!$. However, this notion is not helpful in our present context.
 
For $k>1$, note that $ST^1_{k\ell}$ is the Cayley graph of $Sym_{k\ell}$ generated by the transpositions $(0\;i)$, ($0<i<k\ell)$, but that  
$ST^\ell_k$ is not even a Schreier coset graph of the quotient of $Sym_{k\ell}$ modulo say its subgroup $H_\ell$ generated by the transpositions $(a\;a+1)$, ($0\le a<k$), because the edges of $ST^\ell_k$ are not given by transpositions $(0\;i)$ independently of the values $i$ in different vertices of $ST^\ell_k$. However, Table~\ref{esq} do generalize for every $ST^2_k$, ($k\ge 2$), where the table shows vertically: 

\begin{enumerate}\item[\bf(a)] the right cosets of $V^1_4$ mod the subgroup generated by transpositions $(0\;1),(2\;3)$; \item[\bf(b)] the  representations of such right cosets as vertices of $ST^2_2$; and \item[\bf(c)] assigned generating sets of transpositions $(0\;i)$ per shown right coset of $V^1_3$ or its representing vertex in $ST^2_2$.\end{enumerate}

\begin{table}[htp]
$$\begin{array}{||l||c|c||c|c||c|c||}\hline\hline
Right                &0123 & 2301 &0213 &2031 &0231&2013\\
cosets\; of        &0132 &2310  &0312 &2130 &0321&2103\\
V^1_4         &1023 &3201  &1203 &3021 &1230&3012\\
mod\;H             &1032  &3210 &1302 &3120 &1320&3102 \\\hline\hline
V^2_2         &0011 &1100  &0101 &1010 &0110&1001\\\hline\hline
Gnr.\; set&(0\;2),(0\;3)&(0\;2),(0\;3)&(0\;1),(0\;3)&(0\;1),(0\;3)&(0\;1),(0\;2)&(0\;1),(0\;2)\\\hline\hline
\end{array}$$
\caption{The right cosets of $V^1_4$ as the vertices of $ST^2_2$ and their generating sets.}
\label{esq}\end{table}
 
Tables like Table~\ref{esq}, but for $k>2$, suggest extending the definition of a Schreier coset graph as follows: A {\it Schreier local coset graph} of a group $G$, a subgroup $H$ of $G$ and a generating set $S(Hg)$ for each right coset $Hg$ of $H$ in $G$, is a graph whose vertices are the right cosets $Hg$ and whose edges are of the form $(Hg,Hgs)$, for $g\in G$ and $s\in S(Hg)$. The example in display (\ref{esq}) shows that $ST^2_2$ is a Schreier local coset graph of the group $V^1_4$, its subgroup $H$ generated by the transpositions $(0\;1)$ and $(2\;3)$, and the local generators indicated in the last line of the display. In a similar way, it can be shown for $k>2$ that $ST^2_k$ is a Schreier  local coset graph of $V^2_k$ with respect to its subgroup generated by the transpositions $(2a\;2a+1)$ with $0\le a<k$. Now, the density observed in \cite{D73} must be replaced to be useful in the present context of 2-set star transposition graphs. It is clear that in this sense, the E-sets found in the graphs $ST^2_k$
in Section~\ref{s3} are as dense as they can be, so we say that these E-sets are {\it 2-dense}. Then, the final conclusion of the present section is the following result.

\begin{theorem} 
The E-chain ${\mathcal{ST}}(2)$  of display (\ref{G2}) is  a 2-dense, disjoint segmental, disjoint neighborly E-chain via the E-sets $\Sigma_i^k$ of Theorem~\ref{t2}.
\end{theorem}

\begin{proof} The discussion above in this Section~\ref{Cs} provides all the properties in the statement.
\end{proof}

\section{Pancake 2-set transposition graphs}\label{panqueque}

Let $\pi_i$ be an arbitrary product of independent transpositions on the set $\{1,\ldots,i-1\}$, ($i>1$), where $\pi_1$ and $\pi_2$ are the identity. For each integer $k\ge1$, let $$A(\pi_1,\ldots,\pi_i,\ldots,\pi_{2k-1})=\{(0\;1)\pi_1,\ldots,(0\;i)\pi_i,\ldots,(0\;(2k-1))\pi_{2k-1}\}.$$ Lemma 2 of \cite{D73} implies that for $k\ge 1$ and any choice of the involutions $\pi_i$, ($i\ge 3$), the set $A(\pi_1,\ldots,\pi_{2k-1})$ generates $Sym_{2k-1}$. For each choice of involutions $\pi_1,\pi_2,\ldots$, the sequence of Cayley graphs with generating set $A(\pi_1,\ldots,\pi_{2k-1})$ forms a chain of nested graphs with natural inclusions $\Gamma_k\subset\Gamma_{k+1}$. 

Let $\ell\in\{1,2\}$. If we choose the identity for each entry in $A(\pi_1,\ldots,\pi_{2k-1})$, then we get the $\ell$-set star transposition graphs $ST^\ell_k$. If
$\pi_i=(1\;(i-1))\cdots(\lfloor i/2\rfloor\;\lceil i/2\rceil)$, for $i=3,\ldots,k-1$, then we get the {\it pancake $\ell$-set transpostion graph} $PC_k^\ell$. 
In particular, the pancake 2-set transposition graph $PC^2_k$ has the same vertex set of $ST^2_k$ and its edges involve each the maximal product of concentric disjoint transpositions in any prefix of an endvertex string, including the external transposition being that of an edge of $ST^2_k$.
The graphs $PC_k^1$ were seen in \cite{D73} to form a dense segmental neighborly E-chain 
$\mathcal{PC}(1)=\{PC_1^1,PC_2^1,\ldots,PC_k^1,\ldots\}$. (Figure 2 of \cite{D73} represents the graph $PC^1_4$). In a similar fashion to that of  Section~\ref{Cs}, the following partial extension of that result can be established.
 
\begin{theorem} 
The chain $\mathcal{PC}(2)=\{PC_1^2,PC_2^2,\ldots,PC_k^2,\ldots\}$ is  a 2-dense, disjoint neighborly E-chain via the E-sets $\Sigma_{2k-1}^k$ of Theorem~\ref{t2}, but it fails to be disjoint segmental. A similar result is obtained for any choice of the involutions $\pi_1,\pi_2,\ldots,\pi_i\ldots$ with not all the $\pi_i$s being identity permutations.
\end{theorem}

\begin{proof} Adapting the arguments given for star 2-set transposition graphs in Section~\ref{Cs} can only be done for the E-sets $\Sigma_{2k-1}^k$ in pancake 2-set transposition graphs, since the feasibility for the sets $\Sigma_i^k$, ($1\le i<2k-1$),  to be E-sets is obstructed by the pancake transpositions in $A(\pi_1,\ldots,\pi_{2k-1})$, meaning that we can only establish that the E-chain $\mathcal{PC}(2)$ is dense and disjoint neighborly, but not disjoint segmental. The ``black'" vertices, those whose color is $2k-1$, form an E-set $\Sigma^k_{2k-1}$ with the desired properties, and their removal leaves a $2k-2$-regular graph from which the removal of the ``black" edges, forming an edge subset $E^k_{2k-1}$, leaves the disjoint union of the open neighborhoods $N(v)$ of the vertices $v$ in the E-set $\Sigma^k_{2k-1}$. This behavior is similar for any other choice of the involutions $\pi_1,\pi_2,\ldots,\pi_i\ldots$ with not all the $\pi_i$s being identity permutations, other than $\pi_i=(1\;(i-1))\cdots(\lfloor i/2\rfloor\;\lceil i/2\rceil)$, for $i=3,\ldots,k-1$, which were used precisely to define the pancake graphs.
\end{proof}

\end{document}